\documentclass[11pt]{amsart}
\usepackage{amsmath,amsfonts,amsthm,amssymb,amscd}
\def\classification#1{\def\@class{#1}}
\classification{\null}

\DeclareFontFamily{OT1}{rsfs}{}
\DeclareFontShape{OT1}{rsfs}{n}{it}{<-> rsfs10}{}
\DeclareMathAlphabet{\mathscr}{OT1}{rsfs}{n}{it}

\newcommand{\R}{\mathbb{R}}
\newcommand{\F}{\mathbb{F}}

\newtheorem{lemma}{Lemma}

\newtheorem{theorem}{Theorem}
\theoremstyle{remark}
\newtheorem{remark}[theorem]{Remark}
\newcommand{\oh}{\frac{1}{2}}

\numberwithin{equation}{section}
\title{An explicit incidence theorem in $\mathbb{F}_p$}
\author{Harald Andr\'es Helfgott}
\address{H. A. Helfgott, Department of Mathematics, University of Bristol,
  Bristol BS8 1TW, United Kingdom}
\email{h.andres.helfgott@bristol.ac.uk}
\author{Misha Rudnev}
\address{Misha Rudnev, Department of Mathematics, University of Bristol,
  Bristol BS8 1TW, United Kingdom}
\email{m.rudnev@bristol.ac.uk}
\subjclass[2000]{68R05,11B75}
\begin{document}
\begin{abstract} Let
$P = A\times A \subset \mathbb{F}_p \times \mathbb{F}_p$,
$p$ a prime. Assume that $P= A\times A$ has $n$ elements, $n<p$. See $P$
as a set of points in the plane over $\mathbb{F}_p$. We show that the pairs
of points in $P$ determine $\geq c n^{1 + \frac{1}{267}}$ lines, where $c$
is an absolute constant.

We derive from this an incidence theorem: the number of incidences between a
set of $n$ points and a set of $n$ lines in the projective plane over $\F_p$
($n<\sqrt{p}$) is bounded by $C n^{\frac{3}{2}-\frac{1}{10678}}$, where $C$ is
an absolute constant.
\end{abstract}
\maketitle

\section{Introduction}
$\;\;\;\;$ In 1983 J. Beck proved the following incidence theorem in
$\R^2$. Let $P$ be a set of points in $\mathbb{R}^2$. Then either $P$ contains
$c |P|$ points on a straight line, or the pairs of points of $P$ determine at
least $c |P|^2$ distinct lines, where $c$ is an absolute constant. (We write
$|S|$ for the number of elements of a finite set $S$; $c$ and $C$ will
henceforth denote some absolute constants, changing from one place to
another.) Beck's original paper (\cite{Be}) was followed -- in the same issue
of the same journal --  by a result of Szemer\'edi and Trotter stating that the number of incidences between $m$ straight lines and $n$ points in $\R^2$ is $O\left(m+n+(mn)^{\frac{2}{3}}\right).$ The Szemer\'edi-Trotter theorem implies Beck's theorem as a corollary.

We prove a finite field version of Beck's theorem.

\begin{theorem}\label{delta}
 Let $A\subset \F_p$, where $p$ is larger than an absolute constant. Let
 $P=A\times A \subset \mathbb{F}_p^2$. Let $L = L(P)$ be the set of all lines
 determined by pairs of points of elements of $P$.
If $|A|<\sqrt{p}$, then \[|L(P)|\geq c |P|^{1+\frac{1}{267}},\] where $c>0$ is
an absolute constant.
\end{theorem}
All constants here and later are independent of $p$.

The statement $|L(P)|\geq c |P|^{1+\frac{1}{267}}$ cannot possibly hold in
$\F_p^2$ in full generality,
as $P=\F_p^2$ itself generates only $O(|P|)$ lines. A non-trivial theorem of
Beck type is implicit in the well-known paper of Bourgain, Katz, and
Tao (\cite{BKT}) -- namely, their results imply that, if $P\subset \F_p^2$ is
a Cartesian product and $|P|<p^{2-\Delta}$ ($\Delta>0$), then
$|L(P)|= O\left(|P|^{1+\delta(\Delta)}\right)$ for some
$\delta(\Delta)>0$. However, there was no attempt to make explicit or optimise $\delta(\Delta)$.

The problem of studying $L(P)$ is an incidence problem, and \cite{BKT} proves
that for $n<p^{2-\Delta}$, the number
of incidences between $n$ lines and $n$ points is
$O\left(|n|^{\frac{3}{2}-\delta_1(\Delta)}\right)$.
This statement happens to be on about the same level of generality as claiming
that a
set of $n$ points that is a Cartesian product $A\times A$ (with
$n=|A|^2<p^{2-\Delta}$)
 determines $\Omega\left(n^{1+\delta_2(\Delta)}\right)$ distinct lines.
These results in \cite{BKT} stem from a non-trivial sum-product estimate in
$\F_p$ proven in \cite{BKT}.
 The sum-product estimate of \cite{BKT} (the following formulation also includes Konyagin's contribution \cite{Ko}) says that,
if $A+A$ and $A\cdot A$ denote, respectively, the set of all sums and
products of pairs of elements of $A\subset\F_p$, then, as long as $|A|< p^{1-\Delta},$ one has
\begin{equation}
\max(|A+A|,\;|A\cdot A|)\;\geq c |A|^{1+\delta_3(\Delta)},
\label{spest}\end{equation}
for some absolute $c$. The quantitative relation between the deltas was not
established.

(At this point one traditionally mentions the Erd\"os-Szemer\'edi conjecture,
which states that,
 if $A$ is a subset of integers, then $\max(|A+A|,|A\cdot A|) \geq c |A|^{1+\delta} $
 holds for any $0<\delta<1,$ where $c$ is allowed to depend only on $\delta$, but not on $|A|$.)

We return to $\F_p$ for the rest of the paper. On the level of the existence
of positive exponents, a nontrivial sum-product estimate implies a non-trivial
incidence or
Beck-type theorem, and conversely. Garaev \cite{Ga} succeeded in obtaining a
quantitative sum-product estimate in $\F_p$: for a small enough (say
$|A|<\sqrt{p}$) subset $A$ of $\F_p$, either $A+A$ or $A\cdot A$ has
cardinality at least
$|A|^{1+\frac{1}{14}}$, up to a multiple of a power of $c\log |A|$. Katz and
Shen (\cite{KS}) elaborated on a particular application of the
Pl\"unnecke-Ruzsa inequality  in Garaev's proof and improved the result to
$|A|^{1+\frac{1}{13}}$, up to a multiple of a power of $c \log |A|$. Bourgain
and Garaev (\cite{BG}) incorporated a covering argument (whose variant we cite
as Lemma \ref{cover}) and improved the estimate to $|A|^{1+\frac{1}{12}}$, up to a
multiple of a power of $c \log |A|$. Li (\cite{Li}) showed that a multiple of
a power of $\log |A|$ can be done away with; thus the best result known states
that, for $|A|<\sqrt{p}$,
$$
\max(|A+A|,\;|A\cdot A|)\;\geq \;c |A|^{1+\frac{1}{12}},
$$
where $c$ is an absolute constant.


Our construction uses the techniques laid out in the above-mentioned papers and locally follows rather closely the exposition from some of those papers.
The main result, Theorem 1, an explicit Beck type incidence theorem, implies
the following incidence theorem for $n$ points and $n$ straight lines in the
projective plane
$\mathbb{P}^2(\F_p)$.

\begin{theorem}\label{deltainc}
If $P$ and $L$ are sets of points and lines in $\mathbb{P}^2(\F_p)$ with $|P|, |L| = n < p$, then the number of incidences
$$
I(P,L) = |\{(p,l)\in P\times L:\,p\in l\}| \leq C n^{\frac{3}{2} - \frac{1}{10678}}
$$
for some absolute $C$.
\end{theorem}
We call $\{(p,l)\in P\times L:\,p\in l\}$ the set of {\em incidences}.

The proof of this theorem repeats the pigeonholing argument of \cite{BKT} until it merges with the proof of our Theorem \ref{delta}. It is given at the end of this note.

\begin{remark} In both Theorems \ref{delta} and \ref{deltainc} one can easily
extend the estimate to larger sets, by a straightforward adaptation of Case (i) of the
estimate (\ref{sq}) in the end of our proof.
 We haven't done so aiming at an estimate which does not contain $p$ explicitly;
we leave the case of ``larger'' sets to the interested reader.\end{remark}

\subsection{Acknowledgments}

H. A. Helfgott is supported in part by EPSRC grant EP-E054919/1. The authors would like to thank M. Garaev, T.  Jones and O. Roche-Newton for their helpful comments on the draft of this paper.

\section{Background in arithmetic combinatorics}
We use the following largely standard arithmetic combinatorics lemmata.
In the sequel, in order to suppress constants, we will use
 the $\ll ,\,\gg,\approx$ notations in estimates:
 $|X|\gg|Y|$ means $|X|\geq c|Y|$ for some $c$,
 $|X|\ll|Y|$ means $|X|\leq c|Y|$ for some $c$,
$|X|\approx|Y|$ means $|X|\leq C_1 |Y|$ and $|X|\geq c_2 |Y|$ hold
 for some $C_1$, $c_2$. We abuse the English language in accordance with these notations by saying ``at least'', ``at most'', or approximately in the sense conveyed by the symbols $\gg,\ll,\,\approx$, respectively. To avoid confusion, we will enclose
``at least'' and ``at most'' in quotation marks when we use them in this way.

Thus, for example,
 saying that $|X|$ is ``at least'' $|Y|$ means $|X|\gg|Y|$.

We adopt the following formulation of the Balog-Szemer\'edi-Gowers theorem.
\begin{lemma}[Balog-Szemer\'edi-Gowers theorem]\label{lem:BSG} Let $X,Y$ be additive sets of $n$ elements, and $\alpha\in (0,1)$. Suppose that there is a set of $\alpha n^2$ pairs of elements
$(x,y)\in X\times Y$ on which the sum $x+y$ takes
at most $n$ distinct values.

Then there exist subsets $X'\subseteq X$, $Y'\subseteq Y$, with
$|X'|,|Y'|\gg \alpha n,$ such that
$$
|\text{Range of}\;x+y\mbox{ on } X'\times Y'|\ll  \alpha^{-5} n.
$$\end{lemma}
The modern graph-theoretical proof of the Balog-Szemer\'edi-Gowers
theorem can be found in \cite[Thm.\ 2.29]{TV}. The proof as appears in
that standard reference appears to have a typographical error,
however, which leads to an exponent of $-4$, rather than to the
correct (and somewhat weaker) exponent of $-5$ that we have above.
For a proof yielding the exponent $-5$, see \cite{FS}.

We adopt the following form for the Pl\"unnecke-Ruzsa inequality, due to Ruzsa (\cite{R}).
\begin{lemma}\label{pr}
Let $Y;\; X_1,\ldots X_k$ be additive sets. Then there exists a non-empty
subset $Y'\subseteq Y$, such that
\begin{equation}\label{pr1}
|Y'+ X_1+\ldots+X_k| \;\leq  \;\frac{\prod_{i=1}^k|Y+X_i|}{|Y|^{k-1}}
|Y'|.
\end{equation}
\end{lemma}

Ruzsa's inequality immediately implies that
\begin{equation}\label{pr2}
|X_1+\ldots+X_k| \leq \frac{\prod_{i=1}^k|Y+X_i|}{|Y|^{k-1}},
\end{equation}
for any ``dummy set'' $Y$.
To alow for $-$ signs as well $+$ signs, this inequality is often used in conjunction with the Ruzsa distance inequality
\begin{equation}\label{pr3}
|X_1-X_2| \leq \frac{|X_1-X_3||X_3-X_2|}{|X_3|} .
\end{equation}
(See, e.g., \cite[Lemma 2.6]{TV} for a brief proof of (\ref{pr3}).)

\begin{remark}Katz and Shen (\cite{KS}) showed that, at the expense of acquiring a constant in the left-hand side of (\ref{pr1}), one can make $Y'$ contain an arbitrarily large proportion of $Y$. This trick was used to gain an improvement in (\cite{KS}) and subsequent above-mentioned papers on the sum-product problem. We have not found a way to take advantage of it here, as we shall be dealing with set families, controlling their intersections, refinements therefore being apparently forbidden. I.e., we essentially use Pl\"unnecke-Ruzsa only in the ``crude'' form (\ref{pr2}).\end{remark}

Finally, we need the following covering lemma (see e.g. \cite{BG}, \cite{Sh}, \cite{Li}; the proof is a Cauchy-Schwartz type averaging argument).
\begin{lemma}\label{cover} Let $X_1$ and $X_2$ be additive sets. Then for any $\varepsilon\in (0,1)$ and some constant $C(\varepsilon)$, there exist ${\displaystyle \frac{C(\varepsilon)}{|X_2|}\min(|X_1+X_2|,|X_1-X_2|) }$ translates of $X_2$ whose union contains not less than $(1-\varepsilon)|X_1|$ elements of $X_1$.
\end{lemma}

\section{Proof of Theorem \ref{delta}}
We will prove the result more generally for $P = A_1 \times A_2 \subset
\mathbb{F}_p \times \mathbb{F}_p$, $|A_1| = |A_2| = n$, $n<\sqrt{p}$.

 Let $L(P)$ be the set of straight lines generated by pairs of elements of $P$.
Suppose that
\begin{equation}
|L(P)| \approx n^{2+2\delta}.
\label{ctr}\end{equation}
where $\delta < \frac{1}{267}$. We will show how to reach a contradiction.

The approximately $n^4$ pairs of distinct points of $P$ are distributed between approximately $n^{2+2\delta}$ lines. This implies that a positive proportion of approximately $n^4$ pairs of those points are supported on rich lines, meaning lines
 with ``at least'' $n^{1-\delta}$ points on each. These rich lines
thus contain ``at least'' $n^{5-\delta}$ collinear triples of distinct points of $P$.

Hence, the equation
\begin{equation}
\left|
    \begin{array}{ccc}
      1 & 1 & 1 \\
      x_1 & x_2 & x_3 \\
      y_1 & y_2 & y_3 \\
    \end{array}
  \right| = 0,\qquad x_1,x_2,x_3\in A_1,\; y_1,y_2,y_3\in A_2
  \label{det}\end{equation}
 has ``at least''
  $n^{5-\delta}$  solutions. Then, for some fixed and non-equal
$y_1, y_2 \in A_2$, the equation
  \begin{equation}
\left|
    \begin{array}{ccc}
      1 & 1 & 1 \\
      x_1 & x_2 & x_3 \\
      y_1 & y_2 & y_3 \\
    \end{array}
  \right| = 0,\qquad x_1,x_2,x_3\in A_1,\; y_3\in A_2
   \label{dett}\end{equation}
  has ``at least'' $n^{3-\delta}$  solutions. Since we can translate and dilate
$A_2$, we can assume $y_1=0$, $y_2=1$ without loss of generality.
Then the condition (\ref{dett}) turns into a claim that
 \begin{equation}
  x_1(1-y_3) + x_2 y_3 \in A_1 \qquad x_1,x_2\in A_1,\,y_3\in A_2
  \label{bssetup}\end{equation}
happens for ``at least'' $n^{3-\delta}$ triples $(x_1,x_2,y_3)$. Assuming $y_3\not=0,1$ does not change the situation, as long as $\delta<1$. Let us define $B$ as
the set of elements of the form
  $b=\frac{y_3}{1-y_3}$, $y_3\in A_2\setminus\{0,1\}$. Clearly $|B| \approx |A_2|=n$ and $B$ does not contain zero.

  (The set-up here is similar to that of Theorem C in Bourgain (\cite{B}), which would at this point yield the existence of $\delta$, with a possibility of its quantitative estimate. However, if one chases through Bourgain's arguments, the value of $\delta$ appears to be considerably smaller than what we obtain.)

  Let $B_1$ be the set of popular elements of $B$,
meaning the set of all
  $b\in B$ such the equation (\ref{bssetup}) has ``at least'' $n^{2-\delta}$ solutions $(x_1,x_2,y_3)\in A_1\times A_1\times A_2$ with
$\frac{y_3}{1-y_3} = b$. By the pigeonhole principle,
  $$|B_1|\gg n^{1-\delta}.$$ Besides, fixing $b$ fixes $y_3$, so for each $b\in B_1$ the condition
  $$
  x_1+b x_2 \in \frac{1}{1-y_3(b)} A_1
  $$
  holds for ``at least'' $n^{2-\delta}$ pairs $(x_1,x_2)\in A_1\times A_1$.

  Applying
the Balog-Szemer\'edi-Gowers theorem for each $b$, yields the existence of the subsets  $ A^1_b$ and $A^2_b$ of $A_1$, with $|A^i_b| \gg |A_1|^{1-\delta} = n^{1-\delta},\,i=1,2$, such that
  \begin{equation}
  |A^1_b + bA^2_b|\ll  n^{1+5\delta}.
  \label{bsest}\end{equation}

%
Let us restate that
     \begin{equation}
|A^1_b|,|A^2_b|,|B_1| \gg n^{1-\delta}.\label{sizes}\end{equation}

   Let $A_b = A^1_b\times A^2_b$.  In view of (\ref{sizes}) and by Cauchy-Schwartz, since each $A_b$ is a subset of $A\times A$,
     $$
     |B_1| n^{2-2\delta}\ll  \sum_{b\in B_1}|A_b| \leq n\left(\sum_{b,b'\in B_1} |A_b\cap A_{b'}|\right)^{\frac{1}{2}},
     $$
     so for some $b_*\in B_1$ (which can be assumed non-zero),  denoting $A_* = A_{b_*}$, we have
     \begin{equation}
     \sum_{b\in B_1} |A_b\cap A_{*}| \gg |B_1| n^{2-4\delta}.
     \label{inter}\end{equation}

     Let $B_2$ be a popular subset of $B_1$, namely
a set such that for all $b\in B_2$ we have
     \begin{equation}
     A_{b\wedge *}\equiv |A_b\cap A_{*}|\; \gg \; n^{2-4\delta};
     \label{intersect}
     \end{equation}
     clearly,
     \begin{equation}\label{btwo}
     |B_2| \;\gg \;n^{1-5\delta}.
     \end{equation}
Besides, each $A_{b\wedge *}$ is a cartesian product:
$$
A_{b\wedge *} = A_{b\wedge *}^1\times A_{b\wedge *}^2,
$$
and therefore, for $i=1,2$
     \begin{equation}\label{awedge}
     A_{b\wedge *}^i \equiv  A_{b_*}^i \cap A_b^i\;\gg \;n^{1-4\delta}.
     \end{equation}

We now apply the Pl\"unnecke-Ruzsa inequality
(\ref{pr2})
 with $k=2$, as well as (\ref{bsest}),  and the set cardinality estimates (\ref{sizes}) and (\ref{awedge}) to draw the following conclusions. For each $b\in B_1$ we have
 \begin{equation}\label{sumsetone}\begin{aligned}
   |A^1_{b} + A^1_{b}|  \leq \frac{|A^1_{b} + b A^2_{b}|^2}{|A^2_{b}|} \ll   n^{1+11\delta}, \\
 |A^2_{b} + A^2_{b}|  \leq \frac{|A^1_{b} + b A^2_{b}|^2}{|A^1_{b}|} \ll   n^{1+11\delta}.
\end{aligned}
\end{equation}
Furthermore, for each $b\in B_2$:
$$\begin{aligned}
   |b_* A^2_{b} + b A^2_{b}|  & \leq \frac{|b_* A^2_{b} + b_* A^2_{b\wedge *}
     | |b_* A^2_{b\wedge *} + b A^2_{b}| }{ | A^2_{b\wedge *} |}
\ll n^{15 \delta} |b_* A_{b\wedge *}^2 + b A_b^2|
\\ & \ll
  n^{15\delta}  \frac{|A^1_{b\wedge *} + b_* A^2_{b\wedge *} | |A^1_{b\wedge *} +  b A^2_{b}|}{ | A^1_{b\wedge *} |}\ll n^{1+29\delta},\end{aligned}
 $$
by (\ref{pr2}), (\ref{awedge}), (\ref{sumsetone}), and (\ref{bsest}).

   Then

$$
   |b_* A^2_{*} + b A^2_{b}|  \leq \frac{|b_* A^2_{*} + b_* A^2_{b\wedge *} | |b_* A^2_{b\wedge *} + b A^2_{b}| }{ | A^2_{b\wedge *} |} \ll n^{1+44\delta}.
$$

   and

       \begin{equation}\label{sumsets}
   |b_* A^2_{*} + b A^2_{*}|  \leq \frac{|b_* A^2_{*} + b A^2_{b\wedge *} | |b A^2_{b\wedge *} + b A^2_{b}| }{ | A^2_{b\wedge *} |} \ll n^{1+59\delta}.
   \end{equation}


Throughout the rest of the proof, to save on the number of indices used, let
us refer to  $A^2_*$ as $X$ and to $ b_*^{-1} B_2 $ as $Y$.

Let us use the symbol

\begin{equation}
K =\max_{y\in Y} |X+yX|,\qquad\mbox{so} \qquad K\ll  n^{1+ 59\delta}.\label{ka}\end{equation}

Consider the equation
\begin{equation}\label{pairs}
x_2+ y \tilde x_1 = \tilde x_2 + y x_1 ,\qquad x_1,x_2,\tilde x_1,\tilde x_2\in X,\;y\in Y.
\end{equation}
This equation has at least $\frac{|Y| |X|^4}{K}$ solutions, as follows by applying Cauchy-Schwartz for each  individual $Y$ and then summing over $y\in Y$. Equation (\ref{pairs}) is equivalent to
$$
x_2-\tilde x_2 =  y (x_1 - \tilde x_1),\qquad x_1,x_2,\tilde x_1,\tilde x_2\in X,\;y\in Y.
$$

Hence, for some fixed $(\tilde x_1, \tilde x_2)\in X\times X$, the above equation has at least $\frac{|Y| |X|^2}{K}$ solutions. Let $X_{1} = X- \tilde x_1,\; X_2=X- \tilde x_2$ be translates of $X$ by $\tilde x_1$ and $\tilde x_2,$ respectively. The equation

\begin{equation}\label{lines}
v = y u,\qquad u\in X_1,\;v\in X_2,\;y\in Y
\end{equation}
has at least $\frac{|Y| |X|^2}{K}$ solutions. Consider the set $X_1\times X_2$ in $\F_p^2$. The bound we have just given as to the number of solutions of equation (\ref{lines}) can be rephrased as saying that the set of straight lines through the origin with slopes in $Y$ makes at least $\frac{|Y| |X|^2}{K}$ incidences with $X_1\times X_2$.

In the remainder of the proof will assume that $\frac{|X|^2}{K}\gg 1,$ i.e.,
 that the lines in question contain ``at least'' one point each on average.
This follows immediately from (\ref{sizes}) and (\ref{ka}) once we
assume $\delta< \frac{1}{61}$.

Not less than $50\%$ of the incidences specified by (\ref{lines}) are
 contributed by rich lines with ``at least''
$\frac{|X|^2}{K}$ points thereon. The number of rich lines is not greater than $|Y|$, and ``at least'' $\frac{|Y||X|}{K}$.

Those ``at least'' $\frac{|Y| |X|^2}{K}$ points of $X_1\times X_2$ lying on rich lines can have $|X|$ different abscissae. Hence, there is a vertical
set $u_* \times X_2$
for some non-zero $u_*\in X_1$ intersected by ``at least'' $\frac{|Y| |X|}{K}$ rich lines.

 Thus, we have a subset $Y_1\subset (X_2\cap u_*Y)$ of cardinality \begin{equation}\label{beone}|Y_1|\gg \frac{|Y| |X|}{K}.\end{equation}
In the original notations, $Y_1$ lies in the intersection of $u_* b_*^{-1} B_2$ and some translate of
$A^2_*$; besides,
\begin{equation}
|Y_1| \gg n^{1-65\delta},
\label{e1}\end{equation}
by the bounds (\ref{sizes}), (\ref{btwo}), (\ref{sumsets}).

Let $R$ be the set of all elements expressed via $r= \frac{p-q}{s-t}$, where $p,q,s,t$ are elements of $Y_1$ and $s\neq t$.
Let us consider two cases: (i) $|R|\geq |Y_1|^2$ and
 (ii) $|R|< |Y_1|^2$. Since $n<\sqrt{p}$, $R=\F_p$ is a possibility only in
case (i).

\medskip
Let us consider Case (ii) first.
For any $\xi\not\in R$ and any non-equal $y_1,y_2\in Y_1$, the sum $y_1+\xi
y_2$ has a single realisation, as
\begin{equation}\label{onceagain}y_1 + \xi y_2 = y_1' +\xi y_2',\qquad y_1,y_2,y_1',y_2'\in Y_1\end{equation} would
imply $\xi = \frac{y_1 - y_1'}{y_2' - y_2}$. Thus, for every subset
$Y_1'\subset Y$ and any nonzero $\xi \not \in R$,
\begin{equation}
|Y_1' + \xi Y_1'|=|Y_1'|^2.
\label{sq}\end{equation}

For any  $y\not = 0$ we have some $p,q,s,t\in Y_1$, such that $\xi =
\frac{p-q}{s-t} + y$ lies in the complement of $R$, for otherwise $R+\{y\}=R$,
which is possible only if $R=\F_p$. In particular, this holds for $y=1$.
 Recall that $Y_1$ is a subset of $u_* b_*^{-1} B_2$, and so we may regard  $p,q,s,t$ as elements of $B_2$.
We have, then, some fixed $p,q,s,t\in B_2$ such that
\begin{equation}
|Y_1|^2 \ll |Y_1' + \xi Y_1'|   \leq |Y_1'  + Y_1'+ \frac{p-q}{s-t} Y_1'|,
\label{that}\end{equation}
for any  $Y_1'\subseteq Y_1$, that constitutes a positive proportion of $Y_1$, to be chosen next.

We now use Lemma \ref{cover}.
Let us first show that for any $b\in B_2$ we can cover $99\%$ of the elements
of the sets $bY_1$ (a subset of a translation of $bA^2_*$) or $-bY_1$ (a subset of a translation of $-bA^2_*$) by ``at most'' $n^{24\delta}$ translates of the set
$A^1_*$. Indeed, $A^1_{b\wedge *}=A^1_b\cap A^1_*$ is a subset of $A^1_*$, and by
Lemma \ref{cover} and (\ref{pr2}), we can cover
$99\%$ of the elements of either $bY_1$ or $-b Y_1$ by ``at most''
\begin{equation}
\frac{ |A^1_{b\wedge *}+ bY_1| } {|A^1_{b\wedge *}|} \leq \frac{ |A^1_{b\wedge *}+ bA^2_*| } {|A^1_{b\wedge *}|} \leq\frac{ |A^1_{b\wedge *}+  bA^2_{b\wedge *}||bA^2_{b\wedge *}+  bA^2_{*}| } {|A^1_{b\wedge *}||A^2_{b\wedge *}|} \ll
n^{24\delta}\label{thot}\end{equation}
translates of $A^1_{b\wedge *}$, and hence of $A^1_*$. In the last estimate we've used (\ref{bsest}), (\ref{awedge}), and (\ref{sumsetone}).

This altogether enables us to choose $Y_1'$ as a subset containing at least $50\%$ of $Y_1$, and such that $ (p-q)Y_1'$  gets covered by ``at most'' $n^{48\delta}$ translates of $A^1_*+A^1_*$. Let us now, in the same vein, $\tilde A^2_*$ be a a subset containing at least $50\%$ of $A^2_*$, such that $ (s-t)\tilde A^2_*$  gets covered by ``at most'' $n^{48\delta}$ translates of $A^1_*+A^1_*$. Then we apply Pl\"unnecke-Ruzsa to (\ref{that}) as follows:

\begin{equation}\begin{aligned}
 |Y_1' + Y_1' +  \frac{p-q}{s-t} Y_1' | & \ll
 \frac{|\tilde A^2_* + Y_1'+Y_1'||\tilde A^2_* + \frac{p-q}{s-t} Y_1'|}{ |\tilde A^2_*|} \\ & \ll \frac{|A^2_* +A^2_* +A^2_* |}{n^{1-\delta} } |\tilde A^2_* + \frac{p-q}{s-t} Y_1'| \\ & \ll n^{18\delta}|\tilde A^2_* + \frac{p-q}{s-t} Y_1'|,
\end{aligned}\label{thatt}\end{equation}
after applying Pl\"unnecke-Ruzsa with $k=3$ and the ``dummy set'' $b_*^{-1} A^1_*$, using (\ref{sizes}).

The covering argument above implies that
$$|\tilde A^2_* + \frac{p-q}{s-t} Y_1'| \ll n^{96\delta} | A^1_* + A^1_*+ A^1_*+A^1_*| \ll n^{1+119\delta}, $$
by applying Pl\"nnecke-Ruzsa with the ``dummy set'' $b_*A_*^2$ and $k=4$, using (\ref{bsest}) and (\ref{sizes}).

Therefore returning to (\ref{that}), we have
\begin{equation}
|Y_1'|^2 \ll  n^{1+137\delta}.
\label{e3}\end{equation}
Comparing this with (\ref{e1}), and recalling that
$Y_1'$ contains at least half of the elements of $Y_1$,
 we conclude that  $267\delta\geq 1,$ so $\delta\geq \frac{1}{267}$.
This  ends Case (ii) and essentially ends the proof of Theorem \ref{delta}.

\medskip
Indeed, to analyse Case (i) we observe that if $|R|\geq |Y_1|^2$, then,
summing  the number of ordered quadruples $(y_1,y_2,y_1',y_2')$
of elements of $Y_1$ satisfying equation (\ref{onceagain})  over all
$\xi\in R$ we get no more than
$2|R||Y_1'|^2$ for the total number of solutions, pentuples
$(y_1,y_2,y_1',y_2', \xi)$. Indeed, given $\xi$, the solutions can be either
trivial, with $(y_1,y_2)=(y_1',y_2')$ or not.
The total number of non-trivial solutions is
$|Y_1|^4,$ since such a solution determines $\xi$; the total number of trivial
ones is $|R||Y_1|^2$. Under the assumption of Case (i) the number of trivial
solutions dominates the number of non-trivial ones. Hence, by the pigeonhole
principle and Cauchy-Schwartz, there exists $\xi =\frac{p-q}{s-t}\in R$
such that
\[|Y_1' + \xi Y_1'| \gg \frac{|Y_1'|^4}{|R| |Y_1|^2/|R|} \gg |Y_1'|^2 .\]
for every subset $Y_1' \subset Y_1$ with $|Y_1'|\gg |Y_1|$.

This amounts to a simplification of (\ref{that}), with $\tilde Y'$ replacing
$\tilde Y'+\tilde Y'$ in its right-hand side. The ensuing estimates, done in
exactly the same way as in Case (ii) are therefore better -- by a factor
of $n^{18\delta},$ as follows by inspection of (\ref{thot}) and (\ref{thatt}).

\medskip
This ends the proof of Theorem \ref{delta}. $\Box$

\section{Proof of Theorem \ref{deltainc}}
We follow \cite{BKT}, Section 6.
Suppose, for contradiction, that for some $(P,L)$ one has
\begin{equation}\label{ctr1}
I(P,L)\approx n^{\frac{3}{2}-\epsilon},
\end{equation}
for some $\epsilon>0$. We shall give the lower bound for such $\epsilon$. Note that the role played by the parameter $n$ in this theorem is different from the proof of Theorem \ref{delta}.

Let us first off erase the points in $P$ incident to more than $Cn^{\oh +\epsilon}$ lines of $L$, without changing the notations $(P,L)$. This can be done, as the maximum number of incidences that can come from the set $P_+$ of such points is
$$\begin{aligned}
I(P_+,L) &= \sum_{p\in P_+}\sum_{l\in L} \delta_{pl} \leq \frac{1}{ Cn^{\oh
    +\epsilon}}  \sum_{p\in P_+}\left(\sum_{l\in L} \delta_{pl} \right)^2 \\ &=
\frac{1}{ Cn^{\oh +\epsilon}} \sum_{l,l'\in L} \sum_{p\in P_+}\delta_{pl}\delta_{pl'} \ll  \frac{n^2}{C n^{\frac{1}{2}+\epsilon}},
\end{aligned}$$
since any two distinct lines of $L$ meet at at most a single point of $P_+$.
(Here we use
the notation $\delta_{pl}=1$ if the point $p$ is incident to  the line $l$,
$\delta_{pl}=0$ otherwise.)

This having been done, let $P_1$ be the set of popular points of $P$, each
incident to ``at least'' (recall that saying ``at least'' implies a suitable
constant $c$) $n^{\oh-\epsilon}$ lines of $L$. We have $I(P_1,L)\approx I(P,L)$
by the following argument:
\[I(P,L) - I(P_1,L) = \mathop{\sum_{p\in P}}_{p\notin P_1}
\sum_{\ell \in L} \delta_{p,\ell} < \sum_{p\in P} c n^{\frac{1}{2} - \epsilon}
\leq c n^{\frac{3}{2} - \epsilon} < \frac{1}{2} I(P,L)
\]
for $c$ small enough.

We can now refine $L$ to the subset $L_1$ of popular lines, each incident to
``at least'' $n^{\oh-\epsilon}$ points of $P_1$. By the pigeonhole principle,
it still contributes a positive proportion of incidences. Let $P_2\subset P_1$
with respect to  $L_1$. (That
is,
 $p\in P_1$ is an element of $P_2$ if and only if it is incident to
``at least'' $n^{\oh-\epsilon}$ lines of $L_1$.)

Once again, $I(P_2,L_1)$ satisfies (\ref{ctr1}). (The refinement process could
be iterated any finite number of times, with
the constants obviously getting worse.)


For $p\in P_2$, let $P_p$ be a subset of all points of $P_1$, which are connected to $p$ by some line from $L_1$. By definitions of $P_2,L_1$ we have $|P_p|\gg n^{1-2\epsilon}$, for each $p\in P_2$.
Thus, by Cauchy-Schwartz
$$
|P_2|n^{1-2\epsilon}\ll  \sum_{p\in P_2} |P_p| \leq \sqrt{|P_1|} \sqrt{\sum_{\bar p,\tilde p\in P_2}|P_{\bar p}\cap P_{\tilde p}|},
$$
and so
\[\mathop{\sum_{\bar p, \tilde p \in P_2}}_{\bar p \ne \tilde p}
|P_{\bar p}\cap P_{\tilde p}| = \sum_{\bar p, \tilde p \in P_2}
|P_{\bar p} \cap P_{\tilde p}| - \sum_{\bar p \in P_2} |P_{\bar p}|
\gg \frac{|P_2|^2}{|P_1|} n^{2 - 4 \epsilon} - O(n^2).
\]
Now since $I(P_2,L_1) \approx n^{\frac{3}{2} - \epsilon}$ and each point of
$P_2$
is incident to at most $\ll n^{\frac{1}{2} + \epsilon}$ lines of $L$,
we have $|P_2| \gg n^{1 - 2\epsilon}$, and thus
\[\mathop{\sum_{\bar p, \tilde p \in P_2}}_{\bar p \ne \tilde p}
|P_{\bar p}\cap P_{\tilde p}| \gg \frac{|P_2|^2}{|P_1|} n^{2 - 4 \epsilon}.\]
Therefore one can fix some $(\bar p, \tilde p)\in (P_2\times P_2)$,
$\bar p \ne \tilde p$, such that
$$P_3\equiv |P_{\bar p}\cap P_{\tilde p}|\gg \frac{n^{2-4\epsilon}}{|P_1|} \gg n^{1-4\epsilon}.$$

Each point of $P_3\subseteq P_1$ is incident to ``at least'' $n^{\oh-\epsilon}$ lines of the original set of lines $L$. Without loss of generality, after a projective transformation we can place the points $\bar p$ and $\tilde p$ on the line at infinity, so that the ``at most'' $n^{\oh+\epsilon}$ lines of $L_1$ emanating from these points can be viewed as being parallel to the $x$ and $y$ coordinate axes. In other words, for some $A,B$ of cardinality ``at most''
 $n^{\oh+\epsilon}$ each, we have a subset $P_3$ of $A\times B$ of cardinality
``at least'' $n^{1-4\epsilon}$, such that the number of incidences of $P_3$
with
$L$ is
 \begin{equation}\label{ib}
 I(P_3, L)\gg n^{\frac{3}{2}-5\epsilon}.
 \end{equation}
Now, the number of triples of points of $A\times B$, which are collinear on some line from $L$ is then (by H\"older's inequality or simply noticing that the smallest number of triples is achieved with (\ref{ib}) as equality, with $n$ lines, each supporting the same number of points) ``at least'' $n^{\frac{5}{2}-15\epsilon}\gg |A|^{5-\frac{40\epsilon}{1+2\epsilon}}.$
We can now merge with the proof of Theorem \ref{delta} at its claim (\ref{det}), with $\epsilon = \frac{\delta}{40-2\delta} = \frac{1}{10678}$. $\Box$

\end{document}